\documentclass{amsart}
\usepackage{amssymb,amsmath}
\usepackage{amsfonts}
\usepackage[all]{xy}
\numberwithin{equation}{section}
\begin{document}
\title{ RIEMANN SURFACES AND VERTEX OPERATOR ALGEBRAS}
\author{K.M. Bugajska}
\address{Department of Mathematics and Statistics,
York University,
Toronto, ON, M3J 1P3}
\email{bugajska@yorku.ca}
\date{\today}
\begin{abstract}
We show that for any fixed point $P_0$ on a Riemann surface $\Sigma$ the distinct realizations of cocycles in $\textsl{H}^1(\Sigma,{\mathcal{O}})$ correspond  to the natural appearances of the standard Heisenberg vertex operator algebra $\Pi(P_0)$ and to the commutative Heisenberg vertex operator algebra  $\Pi^0(P_0)$ respectively.
\end{abstract}

\maketitle

\section{Introduction}
The relations between vertex operator algebras (VOA) and Riemann surfaces have a long history  ~\cite{ABK86}, ~\cite{ABM86}. The constructions of vertex algebra bundles over a Riemann surface $\Sigma$ as the associated bundles to the principal $\textit{Aut}\mathcal{O}$-bundle are clearly and beautiful described by Frenkel and Ben-Zvi in ~\cite{EFDB01}. In this note we do not present any new material about vertex operator algebras (VOA) themselves. We merely show how some VOA's occur naturally when we start with a Riemann surface  and how the way of the realization of cocycles $\sigma\in{\textsl{H}^1(\Sigma,{\mathcal{O}})}$ determines the commutative or non-commutative structures of the Heisenberg vertex operator algebras respectively. In other words we show that instead of dealing with the insertion of a vertex operator on a surface $\Sigma$ at $P_0$  a deeper inside into its structure (restricted to an infinitesimal punctured disc centered at a point $P_0$) naturally produces some VOA's .

Any choice of a local coordinate z on $U_0={\Delta_{\epsilon}(P_0)}$ vanishing at $P_0$ determines the decompositions of functions and forms on $U_0$ into their principal and regular parts.  When we have to consider both of this parts (pictures II and III) then the natural dual pairing between $\mathbb{C}((z))dz$ and $\mathbb{C}((z))$ defines antisymmetric product,  or equivalently, a nontrivial cocycle $c(f,g)$ which produces a central extension $(3.3)$ of $\mathbb{C}((z))$. This central extension is the topological Heisenberg  Lie algebra $\mathcal{H}$ and to pass from it to the VOA $\Pi(P_0)$ is simple and standard  ~\cite{EFDB01}.  We simply have to consider the universal enveloping algebra $\textit{U}(\mathcal{H})$, its complition $\widetilde{\textit{U}(\mathcal{H})}$ and notice that $\mathcal{H}$ is an $\widetilde{\textit{U}(\mathcal{H})}$-module. This leads  naturally to the Fock representation $\Pi(P_0)$ of $\mathcal{H}$,  which carries the vertex algebra structure  i.e. all vertex algebra axioms are satisfied. In this point we should notice that the Heisenberg vertex algebra $\Pi(P_0)$ comes with a natural conformal structure with the Virasoro field given by $Y(\frac{1}{2}b_{-1}^2,z)$ and with the central charge $c=1$.  Moreover, the field $Y(b_{-1},z)=b(z)$ corresponds to the mapping of $\mathbb{C}[\Sigma\diagup{P_0}]$ (that occurs in the picture P-II) onto $d\mathbb{C}[\Sigma\diagup{P_0}]$ (that occurs in P-III; see (2)and (3) at the end of the section $(2.4)$) and hence it  provides the equivalence between two meromorphic pictures  P-II and P-III of $\textit{Picc}_0$.

 However, when  (as in the picture I) we realize the cocycles in $\textsl{H}^1(\Sigma,{\mathcal{O}})$ using merely the holomorphic differentials (more exactly the Mittag-Leffler distibutions determined by them) then  only  the subspace $z^{-1}\mathbb{C}[z^{-1}]$ of $\mathcal{K}_{P_0}\cong{\mathbb{C}((z))}$ is involved. Since the cocycle $c(f,g)$ vanishes on this subspace we cannot construct its central extension and we arrive to the commutative Heisenberg vertex operator algebra $\Pi^0$ given by $(3.12)$.
 
 The space of holomorphic functions   $\mathbb{C}((z))$ on $\Delta_{\epsilon}^{*}(P_0)$ is nothing else but the space of the holomorphic sections of the holomorphically trivial line bundle $\underline{\textbf{C}}$ restricted to $\Delta_{\epsilon}^{*}(P_0)$. (Similarly the space $\mathbb{C}((z))dz$  may be identified with the space of sections of the canonocal line bundle restricted to $\Delta_{\epsilon}^{*}(P_0)$.) The decomposition $(2.8)$ into equivalence classes $[\sigma]\in{\textsl{H}^1(\Sigma,{\mathcal{O}})}$ corresponds exactly to the meromorphic picture II of $\textit{Picc}_0$  and the scalar product $(3.2)$ on this space of sections leads to the Heisenberg Lie algebra $\mathcal{H}$.  Now we may ask a question whether there exists a holomorphic line bundle  $\xi\in{\textit{Picc}_n}$, $n>0$, whose holomorphic sections over $\Delta^{*}_{\epsilon}(P_0)$ carry a nice algebraic structure leading to some VOA as in the Heisenberg case. It occurs that when we fix any non-singular even spin bundle $\xi_{\beta}\in{\textit{Picc}_{g-1}}$ then the space $\Gamma(\Delta^{*}_{\epsilon}(P_0),{\mathcal{O}(\xi_{\beta})})$ of holomorphic sections of $\xi_{\beta}$ over $\Delta^{*}_{\epsilon}(P_0)$  (i.e.the space of the half-forms attached at $P_0$) may be identified with the space $H=\oplus{\mathbb{C}\nu_n}$, $\nu_{n}=z^{n}dz^{\frac{1}{2}}$ equipped with the natural symmetric product  ${\left\langle f,g\right\rangle}=Res_{z=0}f(z)g(z)dz$, ${f,g}\in{H}$. This immediately produces the Clifford algebra $\textsl{Cl}H$ with the underlying space isomorphic to the one of the infinite exterior algebra $\Lambda^{\infty}H$.  An irreducible representation of $\textsl{Cl}H$ can be given by a minimal left ideal $\textsl{Cl}H{\Psi_0}$ where ${\Psi_0}=\nu_{-1}\nu_{-2}...$.  This ideal  is naturally isomorphic to the underlying space of the fermionic vertex operator algebra $\Lambda(P_0)$ which itself is isomorphic to the Heisenberg  lattice vertex algebra $V_{\mathbb{Z}}$. 
 
 Similarly as $\mathbb{C}[{\Sigma}\diagup{P_0}]$ provides the sections of the line bundle $\underline{\mathbb{C}}$ with a single pole at $P_0$, the space $W={\mathbb{C}[\Sigma\diagup{P_0}]}\sigma_{\beta}$  consists of all sections of $\xi_{\beta}$ with a single pole at $P_0$ ($\sigma_{\beta}$ denotes a unique section of $\xi_{\beta}$ whose divisor has infinite part equal exactly to $P_{0}^{-1}$). The subspace $W\subset{H}$ corresponds to an  alement $\tau$ of the minimal left ideal $\textsl{Cl}H\Psi_0$ and hence to an element $\tau$ of the component $\Lambda^{(0)}\subset{\Lambda}(P_0)$. Since $\Lambda^{(0)}\cong{\Pi_0}\subset{V_{\mathbb{Z}}}$ the Pluecker relations for $W\in{GrH}$ expressed in the Grassmann variables can be written exactly as the equation $(4.22)$ (in terms of the fields $\Psi(z)$  and $\Psi^{*}(z)$) or equivalently as the equation $(4.23)$ (when $\tau$ is viewed as an element of $\Pi_{0}\subset{V_{\mathbb{Z}}}$).
 
 The Clifford algebra  generated by the sections in $\Gamma({\Delta^{*}_{\epsilon}}(P_0),{\mathcal{O}}(\xi_{\beta}))$ was already introduced in a few other papers (for example ~\cite{AGM86}, ~\cite{BMF87}). Our presentation is only slighty different and is given mainly to illustrate the path ``from a Riemann surface $\Sigma$ to a VOA''. Moreover, since the equations $(4.22)$ and $(4.23)$ produce the whole K-P hierarchy of soliton equations associated to the surface $\Sigma$ (seen as the spectral curve of some differential operator) we describe this problem in the appendix.
 
 This paper is very elementary. We do not even introduce the moduli spaces and the Atiyah algebras of infinitesimal symmetries of vector bundles over them. However, the equivalence of our pictures II and III can be seen as a simple illustration of very sophisticated and famous relation between central charges $c_j=(6j^2-6j+1)c_1$ of the representations of the Virasoro vertex algebra and Chern classes of the determinant line bundle of the vector bundle over moduli space whose fiber over $\Sigma$ is the space of j-differentials on $\Sigma$ $(i.e.  c_0=c_1)$ ~\cite{ABS88}

\section{PRELIMINARIES}

\subsection{Basic Properties of Riemann Surfaces}

Let $\textbf{M}$ be an oriented compact two dimentional smooth real manifold of genus $g\geq{0}$. When $\textbf{M}$ is given a complex structure then it becomes a Riemann surface, say, $\Sigma$. The short exact sequence of sheaves over $\Sigma$ given by
\begin{equation}
0\rightarrow{\mathbb{C}}\stackrel{i}{\rightarrow}{\mathcal{O}}\stackrel{\partial}{\rightarrow}{\mathcal{O}^{1,0}}\rightarrow{0}
\end{equation}

produces the following exact sequence
\begin{equation}
0\rightarrow{\Gamma(\Sigma,{\mathcal{O}^{1,0}})}\stackrel{\delta}{\rightarrow}{\textsl{H}^1(\Sigma,{\mathbb{C}})}\stackrel{i^{*}}{\rightarrow}{\textsl{H}^1(\Sigma,{\mathcal{O}})}\rightarrow{0}
\end{equation}

It  reveals the decomposition of the topological invariant $\textit{H}^1(\textbf{M},{\mathbb{C}})\cong{\textit{H}^1(\Sigma,{\mathbb{C}})}$ into two $g$-dimentional subspaces.  One of them contains   cocycles $[\sigma]={\delta(\omega)}$  that come naturally from  holomorphic differentials  and produce  flat bundle representatives of the holomorphically trivial line bundle $\underline{\mathbb{C}}$ over $\Sigma$  (equivalently  $i^{*}[\sigma]=0$ in $\textsl{H}^1(\Sigma,{\mathcal{O}})$).  The other consists of cocycles with $i^{*}[\sigma]\neq{0}$  which determine  flat bundle realizations of a holomorphically nontrivial line bundles (but still topologically trivial). The set of all classes of holomorphic line bundles (HLB) $\xi$'s over $\Sigma$ is in the one-one correspondence with the set of elements of  $\textsl{H}^1(\Sigma,{\mathcal{O}^{*}})$ and the exact sequence
\begin{equation}
0\rightarrow{\textsl{H}^1(\Sigma,{\mathbb{Z}})}\stackrel{i^{*}}{\rightarrow}{\textsl{H}^1(\Sigma,{\mathcal{O}})}\stackrel{e^{*}}{\rightarrow}{\textsl{H}^1(\Sigma,{\mathcal{O}^{*}})}\stackrel{c}{\rightarrow}{\mathbb{Z}}\rightarrow{0}
\end{equation}

tells us that all HLB's which admit a flat realizations come from cocycles in $\textsl{H}^1(\Sigma,{\mathcal{O}})$  (here $c(\xi)$ denotes the Chern class $c_1(\xi)\in{\mathbb{Z}}$). From the sequences $(2.2)$ and $(2.3)$ we see that the set of such HLB's can be identified with the double quotient 
\begin{equation}
\delta\Gamma(\Sigma,{\mathcal{O}^{1,0}})\diagdown{\textsl{H}^1(\Sigma,{\mathbb{C}})}\diagup{\textsl{H}^1(\Sigma,{\mathbb{Z}})}
\end{equation}
which is equivalent to
\begin{equation}
\textsl{H}^1(\Sigma,{\mathcal{O}})\diagup{\textsl{H}^1(\Sigma,{\mathbb{Z}})}=:\textit{Picc}_0
\end {equation}

Let $\textit{Div}\Sigma$ denote the group of divisors on $\Sigma$. Let $\textit{Div}_0{\Sigma}$ be its subgroup of degree zero divisors and let $\textit{Div}{\mathcal{P}}$ denote the principal divisors (i.e. coming from meromorphic functions). If $\mathcal{M}$ denote the sheaf of germs of meromorphic functions on $\Sigma$ then the vanishing of $\textsl{H}^1(\Sigma,{\mathcal{M}})$ for any compact Rieman surface implies that any HLB $\xi$ in $\textsl{H}^1(\Sigma,{\mathcal{O}^{*}})$ admits a meromorphic section as well as it leads to the following exact sequence
\begin{equation}
0\rightarrow{\Gamma(\Sigma,{\mathcal{O}^{*}})}\rightarrow{\Gamma(\Sigma,{\mathcal{M}})}\stackrel{div}{\rightarrow}{\textit{Div}\Sigma}\rightarrow{\textsl{H}^1(\Sigma,{\mathcal{O}^{*}})}\rightarrow{0}
\end{equation}

Hence, any bundle $\xi$ corresponds to a unique class of divisors, i.e. we have $\textsl{H}^1(\Sigma,{\mathcal{O}^{*}})\cong{\textit{Div}\Sigma}\diagup\textit{Div}{\mathcal{P}}$ and  $\textit{Div}_0{\Sigma}\diagup{\textit{Div}{\mathcal{P}}}$  is  equivalent to $\textit{Picc}_0$.

Summarizing, by fixing a concrete complex structure on $\textbf{M}$ we are splitting the topological invariant $\textsl{H}^1(\textbf{M},{\mathbb{C}})\cong{\textsl{H}^1(\Sigma,{\mathbb{C}})}$ into two types of elements. Cocycles $[\sigma]$'s in $\textsl{H}^1(\Sigma,{\mathbb{C}})$ that produce flat representations of the holomorphic line  bundle $\underline{\mathbb{C}}$ (and naturally  determined by holomorphic differentials)  and cocycles $[\sigma]$'s that represents holomorphically nontrivial line bundles in $\textit{Picc}_0$. The former ones  occur in exactly the same way in all  three pictures below and hence we will not bother with them. However, cocycles whose image $i^{*}[\sigma]\neq{0}$ admit a few, completely different from each other, realizations.

\subsection{Holomorphic Picture}
In this picture we will express all elements $[\sigma]\in{\textit{H}_1(\Sigma,{\mathbb{C}})}$ and hence, all $i^{*}[\sigma]\in{\textit{H}_1(\Sigma,{\mathcal{O}})}$, using only holomorphic differentials. Let ${\omega}\in{\Gamma(\Sigma,{\mathcal{O}^{1,0}})}$. On the one side it determines  the mentioned earlier $\delta(\omega)$ but on he other side    it   produces a representative of a nontrivial element in
\begin{equation*}
\textit{H}_1(\Sigma,{\mathbb{C}})\diagup\delta{\Gamma(\Sigma,{\mathcal{O}^{1,0}})}\cong{\textit{H}_1(\Sigma,{\mathcal{O}})}
\end{equation*} 

as follows. Let us fix an arbitrary point $P_0$ on $\Sigma$ and let $\mathcal{U}=\{U_0,U_1\}$ be an open covering with $U_0$ given by a disc $\Delta_{\epsilon}(P_0)$ centered at $P_0$ and with ${U_1}={\Sigma{\backslash}P_0}$. Let z be a local coordinate on $U_0$ vanishing at $P_0$ and let a holomorphic differential $\omega$ be given by $g(z)dz$ on $U_0$ (where $g(z)$ is a holomorphic function without zeros on ${U_0}\cap{U_1}=\Delta_{\epsilon}^{*}$). This allows us to introduce a Mittag-Leffler (M-L) distribution ${\eta}=(\eta_0,\eta_1)$  subordinate to the covering $\{U_0,U_1\}$ by
\begin{equation*}
{\eta}_0(z)=\frac{1}{zg(z)} \quad z\in{U_0}, \qquad {\eta_1}=const \quad on \quad U_1
\end{equation*}

 This distribution is without solution and it defines  a unique element of $\textit{H}^1(\Sigma,{\mathcal{O}})$. Hence,  all flat bundles representing $\underline{\mathbb{C}}$  are given by the periods of holomorphic differentials and holomorphically nontrivial line bundles which possess flat representatives are given by (M-L distributios determined by) the holomorphic differentials as well. Equivalently we may say that the holomorphic picture is tied to the following short exact sequence
\begin{equation}
0\rightarrow\Gamma(\Sigma,{\mathcal{O}^{1,0}})\rightarrow{\textit{H}^1(\Sigma,{\mathbb{C}})}\rightarrow{\Gamma(\Sigma,{\mathcal{O}^{1,0}})^{*}}\rightarrow{0}
\end{equation}

We will see that this picture (denoted as I) leads to the commutative Heisenberg  vertex operator algebra.

\subsection{Function Picture}
In this picture we realize elements $[\sigma]\in{\textit{H}^1(\Sigma,{\mathcal{O}})}$ using the standard definition of cohomology. We use the same open covering $\mathcal{U}=\{U_0,U_1\}$  and local coordinate z on $U_0$ as before. Now, any cocycle ${\sigma}\in{\mathcal{Z}^1(\mathcal{U},{\mathcal{O}})}$ is given by some holomorphic function  $\sigma(z)\in{\Gamma({U_0}\cap{U_1},{\mathcal{O}})}\cong{\mathcal{K}_{P_0}}$ and its cohomology class $[\sigma]$ corresponds to 
\begin{equation}
[\sigma]=\{\sigma+f+g \quad | f\in{\Gamma(U_0,{\mathcal{O}})}\cong{\mathcal{O}_{P_0}},\quad  g\in{\Gamma(U_1,{\mathcal{O}})}\cong{\mathbb{C}[{\Sigma}\backslash{P_0}]}\}
\end{equation}

Hence $\textit{H}^1(\Sigma,{\mathcal{O}})\cong{\mathbb{C}[\Sigma\backslash{P_0}]}\diagdown{\mathcal{K}_{P_0}}\diagup{\mathcal{O}_{P_0}}$.

\subsection{Form Picture}
Let  $\mathcal{M}^{(1)}_0$ denote the sheaf of germs of meromorphic differentials that have residue equal to zero at each point of $\Sigma$. The exact sequence of sheaves $0\rightarrow{\mathbb{C}}\stackrel{i}{\rightarrow}{\mathcal{M}}\stackrel{d}{\rightarrow}{\mathcal{M}^{(1)}_0}\rightarrow{0}$ produces the exact sequence
\begin{equation}
0\rightarrow{\Gamma(\Sigma,{\mathbb{C}})}\rightarrow{\Gamma(\Sigma,{\mathcal{M}})}\rightarrow{\Gamma(\Sigma,{\mathcal{M}^{(1)}_0})}\rightarrow{\textit{H}^1(\Sigma,{\mathbb{C}})}\rightarrow{0}
\end{equation}

which shows that any class $[\sigma]\in{\textit{H}^1(\Sigma,{\mathbb{C}})}$ comes from a meromorphic differential, say $\mu$, without residue.  When $\mu$ is holomorphic, i.e. ${\mu}={\omega}$, then $[\sigma]={\delta(\omega)}$ forms the trivial cocycle in $\textit{H}^1(\Sigma,{\mathcal{O}})$.  When ${\mu}=df$, $f\in{\Gamma(\Sigma,{\mathcal{M}})}$ then $[\sigma]=[\{\sigma_{\alpha\beta}=0\}]$. Hence 
\begin{equation}
\textit{H}^1(\Sigma,{\mathcal{O}})\cong{\frac{\mathcal{M}^{(1)}_{0}-\Gamma(\Sigma,{\mathcal{O}^{1,0}})}{d\Gamma(\Sigma,{\mathcal{M}})}}
\end{equation}

Let $L(P_0^{-2g})\subset{\mathbb{C}[\Sigma\backslash{P_0}]}$ be the  $g+1$-dimentional vector space of functions whose single pole at $P_0$ has order $\leq{2g}$, that is,  $L(P_0^{-2g})={\mathbb{C}+Span_{\mathbb{C}}\{f_1,\ldots,f_g\}}$  and $f_i$'s are non-gap functions. Let $\Omega_0(P_0^{-2g-1})$  denote the vector space of meromorphic differentials on $\Sigma$ with at most single pole at $P_0$ of order $\leq{2g+1}$
\begin{equation*}
\Omega_0(P_{0}^{-2g-1})={\Gamma(\Sigma,{\mathcal{O}}^{1,0})}\oplus{Span_{\mathbb{C}}}\{\tau^{(2)}_{P_0},\ldots,\tau^{2g+1}_{P_0}\}
\end{equation*}

where the meromorphic differentials $\tau^{(n)}_{P_0}$ have  all a-periods zero and the principal part at $P_0$ equal to $z^{-n}dz$. Since any meromorphic differential with a single pole at $P_0$ of order $>2g+1$ can be written as $\mu+df$ for some ${\mu}\in{\Omega_{0}}(P_0^{-2g-1})$ and for some $f\in{\mathbb{C}}[\Sigma\backslash{P_0}]$ we obtain
\begin{equation*}
 \textit{H}^1(\Sigma,{\mathbb{C}})\cong{\frac{\Omega_0(P^{-2g-1})}{dL(P^{-2g}_0)}}\cong{\Gamma(\Sigma,{\mathcal{O}}^{1,0})}\oplus{\frac{\mathcal{M}^{(1)}(P_0^{-2g-1})}{dL(P^{-2g}_0)}}
 \end{equation*}
 
 and hence 
 \begin{equation}
 \textit{H}^1(\Sigma,{\mathcal{O}})\cong{\frac{\mathcal{M}^{(1)}(P_0^{-2g-1})}{dL(P_0^{-2g})}}\cong{\frac{\mathcal{M}^{(1)}_{P_0}}{d\mathbb{C}[\Sigma\backslash{P_0}]}}
 \end{equation}
 
 (Here $\mathcal{M}^{(1)}_{P_0}$ and $\mathcal{M}^{(1)}(P_0^{-2g-1})$ denote pure meromorphic differentials with a single pole at $P_0$.)
 
   Our local coordinate z at $P_0$ determines the   isomorphisms:  $\mathcal{K}_{P_0}\cong{\mathbb{C}}((z))$,     $\Omega_{\mathcal{K}}\cong{\mathbb{C}}((z))dz$ and  $\Omega_{\mathcal{O}}\cong{\mathbb{C}}[[z]]dz$ and hence allows us to see   the following, mutually equivalent, three basic realizations as
 \begin{enumerate}
 \item  ${\textit{H}^1(\Sigma,{\mathcal{O}})}\cong{\Gamma(\Sigma,{\mathcal{O}^{1,0}})^{*}}$  (in the holomorphic picture P-I)
 \item ${\textit{H}^(\Sigma,{\mathcal{O}})}\cong{\mathbb{C}[\Sigma\backslash{P_0}]\diagdown{\mathcal{K}_{P_0}}\diagup{\mathcal{O}_{P_0}}\cong{\mathbb{C}[\Sigma\backslash{P_0}]}\diagdown{\mathbb{C}((z))}\diagup{\mathbb{C}[[z]]}}$  (in the function picture P-II)
 \item ${\textit{H}^1(\Sigma,{\mathcal{O}})}\cong{d\mathbb{C}[\Sigma\backslash{P_0}]\diagdown{\Omega_{\mathcal{K}}}\diagup{\Omega_{\mathcal{O}}}\cong{d\mathbb{C}[\Sigma\backslash{P_0}]\diagdown{\mathbb{C}}((z))dz}\diagup{\mathbb{C}[[z]]dz}}$ (in the form picture P-III).
 \end{enumerate}
 The  pictures II and III are called ``meromorphic''  since  we have to work with both, the principal parts and  the regular parts of functions or differentials respectively. 

\section{FROM PICTURES TO VOA'S}
\subsection{Noncommutative Heisenberg vertex algebra}
Let $\{1=n_1<n_{2}<...<n_{g}<2g\}$ denote the gap sequence at $P_{0}\in{\Sigma}$  ~\cite{FIK92} .
The realization $(2)$ of $\textit{H}^{1}(\Sigma,{\mathcal{O}})$ tells us that we may identify the double quotient ${\mathbb{C}[\Sigma\backslash{P_0}]\diagdown{\mathbb{C}((z))}\diagup{\mathbb{C}[[z]]}}$  with the space $Span_{\mathbb{C}}\{z^{-n_i}|i=1,\ldots,g\}$.
Similarly, the realization $(3)$ means that  $d{\mathbb{C}[\Sigma\backslash{P_0}]}\diagdown{\mathbb{C}((z))dz}\diagup{\mathbb{C}[[z]]dz}$ can be identified with the space  ${{Span}_{\mathbb{C}}\{z^{-n-1}dz\}}$, or equivalently, with the space ${Span_{\mathbb{C}}}\{\tau^{(n+1)}_{P_0}\}$, where again  $n=n_i$ must belong to the  gap sequence at ${P_0}$. So, the equivalence of P-II and P-III is tied to the  natural mapping d from ${\mathcal{K}_{P_0}}$ to ${\Omega_{\mathcal{K}}}$
\begin{equation*}
  f(z)\stackrel{d}{\rightarrow}{\mu^{f}(z)=f'(z)dz}, \quad f(z)\in{\mathcal{K}_{P_0}} 
\end{equation*}

Using the fact that $\mathbb{C}((z))dz$  is dual to $\mathbb{C}((z))$ with the dual pairing given by the residue  
\begin{equation}
{\left\langle \mu,f\right\rangle}=Res_{z=0}f\mu, \qquad {\mu}\in{\Omega_{\mathcal{K}}}, \quad f\in{\mathcal{K}}_{P_0}
\end{equation}

we may introduce the natural product on $\mathcal{K}_{P_0}$ given by
\begin{equation}
(f,g):={\left\langle \mu^{f},g\right\rangle}=Res_{z=0}gdf=-Res_{z=0}fdg \qquad and \quad hence\quad (z^n,z^m)=n\delta_{n,-m}
\end{equation}

Since this product is antisymmetric it forms a cocycle $c(f,g)=-Res_{z=0}fdg$ which defines a central extension of $\mathbb{C}((z))$:
\begin{equation}
0\rightarrow{\mathbb{C}\textbf{1}}\rightarrow{\mathcal{H}}\stackrel{\pi}{\rightarrow}{\mathbb{C}((z))}\rightarrow{0}
\end{equation}

This extension $\mathcal{H}$ is called the Heisenberg (or oscillator) algebra. It is a complete topological Lie algebra with topological generators $b_n; n\in{\mathbb{Z}}$ and $\textbf{1}$ whose images under the projection $\pi$ above are
\begin{equation}
\pi(b_n)=z^n, \qquad  \pi(\textbf{1})=0
\end{equation}

The cocycle $c(z^n,z^m)=n\delta_{n,-m}$ determines the following commutation relations in $\mathcal{H}$:
\begin{equation}
[b_n,b_m]=n\delta_{n,-m}\textbf{1}, \qquad [\textbf{1},b_n]=0, \quad n,m\in{\mathbb{Z}}
\end{equation}

For any $f=\sum_{k=-K}^{N}{f_{k}b_k}\in{\mathcal{H}}$ we have $\pi(f)=\sum_{k=-K}^{N}{f_{k}z^k}\in{\mathbb{C}[z^{-1},z]}\subset{\mathbb{C}((z))}$. 

  In the commutative algebra ${\mathbb{C}((z))\cong{\pi(\mathcal{H})}}$ the multiplication by $z^{-1}$ of $z^k$ , $k\in{\mathbb{Z}}$, differs only by  a multiplicative constant k from the operation $\partial_{z}z^k$.  So, if we consider $\mathbb{C}((z))$ as a set with the binary operations given by the addition then we may view $z^{-1}$ as an element of $End\mathbb{C}((z))$ whose action, denoted by $''\circ''$, is given by

\begin{equation}
z^{-1}\circ{z^k}=kz^{-1}z^k, \quad k\in{\mathbb{Z}}, \quad so \quad that \quad z^{-1}\circ{f(z)}=f'(z)
\end{equation} 

 Analogously, for any $n>0$ we define the $''\circ''$-action of $z^{-n}$ on $\mathbb{C}((z))$ given by infinite many products of $z^{-n}$ with $z^k$, $k\in{\mathbb{Z}}$ with multipliers depending on k so that the action $z^{-n}\circ{f(z)}=f^{(n)}(z)$ produces the n-th derivative of $f(z)$. In this way, for each $n>0$ we have introduced the mapping 
\begin{equation}
{\kappa}:z^{-n}\rightarrow{z^{-n}{\circ}}\in{End{\mathbb{C}}((z))}
\end{equation} 

Hence $\kappa$ is a one-to-one mapping from $z^{-1}\mathbb{C}[z^{-1}]$ into a  commutative subalgebra of  endomorphisms of $\mathbb{C}((z))$.  Notice that for $n>0$ we cannot associate to $z^{n}={\pi}(b_n)$ any concrete endomorphism  of $\mathbb{C}((z))$  (the integration operation involves arbitrary constants of integration).

 Let us return to the non-commutative Heisenberg algebra $\mathcal{H}$.  Each element ${b_n}\in{\mathcal{H}}$  may be naturally  viewed as an element of $End{\mathcal{H}}$  with the action $b_{n}(f)$ (for  $f={\sum}{f_{k}b_k}\in{\mathcal{H}}$)  given by the commutator i.e. $b_n(f)=[b_n,f]={\sum}{f_k}{\delta}_{k,-n}\textbf{1}$.  Let us introduce an element $I\in{End{\mathcal{H}}}$ by the requirement that the following diagram is commutative
\[
\xymatrix{
b_n\ar[d]_{I} \ar[r]^{\pi}
&z^{n}\ar[d]^{{z^{-1}}{\circ}}\\
[I,{b_n}]\ar[r]^{\pi}
&{nz^{-1}}{z^n}\\
}
\]

This  simply means that we must have $[I,{b_n}]=nb_{n-1}$. Usually we work with the endomorphism $T=-I$ instead. So, for $f={\sum}{f_k}{b_k}\in{\mathcal{H}}$ we have
\begin{equation}
[f,T]={\sum}{f_k}k{b_{k-1}}\stackrel{\pi}{\rightarrow}{f'(z)}\quad equivalently \quad [T,{b_k}]=-kb_{k-1}
\end{equation} 

Let $\pi_{-}$ denote the restriction of the projection $\pi$ in $(3.3)$ to the negative part $\mathcal{H}_{-}$ of $\mathcal{H}$. Let $\pi_{-}'$ denote the composition of $\kappa$ and ${\pi_{-}}$:
\begin{equation*}
b_{-n}\stackrel{\pi}{\rightarrow}z^{-n}\stackrel{\kappa}{\rightarrow}{z^{-n}\circ}=\pi_{-}'(b_{-n})\in{End{\mathbb{C}((z))}} 
\end{equation*}

In particular $\pi_{-}'(b_{-1})={z^{-1}\circ}$  and it maps  $f(z)\rightarrow{f'(z)}$. Now it is natural to define an element $Y(b_{-1},z)\in{End{\mathcal{H}[[z^{\pm}]]}}$ by the property that for each $f={\sum}f_{k}b_{k}\in{\mathcal{H}}$  the action of $Y(b_{-1},z)(f)$  given by the commutator i.e. $Y(b_{-1},z)(f):=[f,Y(b_{-1},z)]$  is as follows:
\begin{equation}
[f,Y(b_{-1},z)]=f'(z)\textbf{1} 
\end{equation}

where $f(z)=\pi(f)$ is in a general case a distribution. This requirement uniquely determines $Y(b_{-1},z)$  as  ${\sum}_{k\in{\mathbb{Z}}}b_{k}z^{-k-1}$.  Thus, for any $f(z)=\pi(f)$ with $f={\sum}f_{n}b_{n}$ we have
\begin{equation}
\xymatrix{
b_{-1} \ar[d]_{Y} \ar[r]^{\pi}
&z^{-1}\ar[r]^{\kappa}
&z^{-1}{\circ}:f(z)\ar[r]
&f'(z)\\
Y(b_{-1},z):f\ar[r]
&f'(z)\textbf{1}\\
}
\end{equation}

Analogously, for $n>0$  we determine $Y(b_{-n},z)\in{End{\mathcal{H}}[[z^{\pm}]]}$  by the correspondence 
\begin{equation}
\xymatrix{
b_{-n}\ar[d]_{Y} \ar[r]^{\pi}
&z^{-n}\ar[r]^{\kappa}
&z^{-n}\circ:f(z)\ar[r]
&f^{(n)}(z)\\
(n-1)!Y(b_{-n},z):f\ar[r]
&f^{(n)}(z)\textbf{1}\\
}
\end{equation}

Usually we denote $Y(b_{-1},z)$ by $b(z)$  which implies that, for $n>1$ an operator $Y(b_{-n},z)$  can be formally written as $\frac{1}{(n-1)!}\partial^{(n-1)}b(z)$.  

 The relations $(3.3)-(3.9)$ naturally lead to the introduction of the Fock representation of the Heisenberg algebra $\mathcal{H}$ in  $\Pi(P_0)={\mathbb{C}[b_{-1},b_{-2},\ldots]}|0>$ together with its natural vertex operator algebra (VOA) structure  ~\cite{VKA97}  . (Vector $|0>$  is killed by the half of the Heisenberg algebra i.e. $b_n|0>=0$ for $n\geq{0}$ and $\textbf{1}|0>=|0>$ and  hence  ${\Pi(P_0)}\cong{\mathbb{C}[b_{-1},b_{-2},\ldots]}$) In this way, for any fix point $P_{0}\in{\Sigma}$ and for any local coordinate $z$ on $U_{0}=\Delta_{\epsilon}(P_0)$ with $z(P_0)=0$ picture P-II ( or the equivalence of pictures  P-II and P-III) results in attaching a VOA $\Pi(P_0)$  at the point $P_0$. Since $\Pi(P_0)$ is produced merely by the Laurent expansions of functions and forms at $P_0$ together with their dual pairing $(3.1)$ this structure is exactly the same for any choice of point $P\in{\Sigma}$ as well as for any other compact Riemann surface.  To recover our surface $\Sigma$ out of $\Pi(P_0)\cong{\Pi}$ we must introduce the conformal blocks ~\cite{YZH93}   $\mathcal{C}(\Sigma,P_{0},z,\Pi)$ which are contained in the restricted dual $\Pi^{*}(P_0)\cong{\Pi^{*}}$ and which are determined by the spaces $\Omega_0(P_{0}^{-2g-1})$ and $L(P_{0}^{-2g})$, introduced in the previous section, respectively. 
 
 \subsection{VOA associated to the holomorphic picture}
 Let point $P_{0}\in{\Sigma}$ and let a local coordinate z be the same as in the previous subsection. Let $\{\omega_1,\ldots,\omega_g\}$ be a basis for the space of holomorphic differentials on $\Sigma$ that is dual to a homology basis $\{a_1,\ldots,b_g\}$. Let $\omega_i(z)={\sum}_{k=0}^{\infty}{\alpha_k^{i}z_{k}dz}$ be the expansion of $\omega_i$ at $P_0$. Let us construct a $g\times{\infty}$ matrix $\textbf{B}^{\infty}$ whose i-th row is given by the coefficients of $\omega_i$  i.e. it is equal to $(\alpha^i_0,\alpha^i_1,\alpha^i_2,\ldots)$. Let $\textbf{B}^{\infty}_N$ denote the matrix obtained from the first N columns of $\textbf{B}^{\infty}$. 
 
 We see immediately that the null space of the matrix $\textbf{B}^{\infty}_{2g}$ is spanned by the principal parts of the non-gap functions for $P_0$. When $\{n_{1}=1<n_{2}<\ldots<n_{g}\leq{2g-1}\}$ is the gap sequence at $P_0$ then for each $n_i$ there exists a holomorphic differential ${\omega}\in{\Gamma(\Sigma,{\mathcal{O}}^{1,0})}$ with zero of degree $n_{i}-1$ at $P_0$  i.e. with the expansion $\omega(z)={\sum}_{k=n_{i}-1}^{\infty}\alpha_{k}z^{k}dz$  at this point. The   Mittag-Leffler  distribution corresponding to $\omega$, is totally determined by $\eta_0(z)$ on $U_0$ which is given by
 \begin{equation*}
 \eta_{0}(z)= g_{-n_{i}+1}z^{-n_{i}}+g_{-n_{i}+2}z^{-n_{i}+1}+\ldots+g_{0}z^{-1}+ reg
 \end{equation*} 
  
  If $\underline{\textbf{d}}(\omega)$ denote a vector in $\mathbb{C}^{2g}$ with components $(g_0,g_{-1},\ldots,g_{-n_{i}+1},0,\ldots,0)$ then the matrix $\textbf{B}^{\infty}_{2g}$ determines the splitting (corresponding to $(2.2)$ and $(2.7)$) of the space $\mathbb{C}^{2g}$ into the g-dimentional subspace spanned by the principal parts of the non-gap (at $P_0$) functions and the g-dimentional subspace spanned by the vectors $\underline{\textbf{d}}(\omega)$. The null space of the whole matrix $\textbf{B}^{\infty}$ is given by the principal parts of all functions ${f}\in{\mathbb{C}[\Sigma\backslash{P}_0]}$.
  
  We see that when we work  with holomorphic differentials exclusively (i.e. with the holomorphic picture P-I) then we work only with principal parts of functions and with the vectors $\underline{\textbf{d}}(\omega)$ for $\omega\in{\Gamma(\Sigma,{\mathcal{O}^{1,0}})}$. Since the cocycle (determined by the residue pairing and equivalent to the antisymmetric product on  $\mathbb{C}((z))$) vanishes on this subspace we cannot construct its nontrivial extension analogous to $(2.3)$ . However, this commutative algebra naturally carries the structure of the commutative Heisenberg vertex operator algebra $\Pi^0$ i.e.
  \begin{equation}
  {\Pi^0(P_0)}\cong{Sym{\mathcal{K}_{P_0}}}\diagup{\mathcal{O}_{P_0}}\cong{Sym{\mathbb{C}((z))}\diagup{\mathbb{C}[[z]]}}
  \end{equation}

  In other words, the commutative vertex operator algebra $\Pi^0(P_0)$ is associated to the holomorphic picture P-I in which any cocycle $[\sigma]\in{\textit{H}^1(\Sigma,{\mathcal{O}})}$ is produced out of a holomorphic differential $\omega$ using the M-L distribution  $\eta(\omega)$. This implies that the space of coinvariants that rediscovers $\Sigma$ from ${\Pi^0(P_0)}\cong{\Pi^0}\cong{Sym({\mathbb{C}((z))\diagup{\mathbb{C}}[[z]]})}$ can be viewed, ~\cite{EFDB01}, ~\cite{TEG06}, as the ring $Fun(\Gamma(\Sigma,{\mathcal{O}^{1,0}}))$ of polynomial functions on $\Gamma(\Sigma,{\mathcal{O}^{1,0}})$.
  
\section{Riemann surface $\Sigma$ and a lattice VOA}
\subsection{Algebras $V_{\mathbb{Z}}$ and $\Lambda$ }
 Suppose that instead of working merely with elements  of $Picc_0$  we would like to include into our considerations all classes of holomorphically equivalent line bundles over $\Sigma$.  Since in this case  we must also  work with the constant sheaf $\underline{\mathbb{Z}}$ over $\Sigma$ (see the sequence $(1.3)$) it could  suggest that we should consider the  Heisenberg VOA's associated to the integral lattice   $\mathbb{Z}$. However this naive reasoning is quite misleading and it is completely false. To understand this problem better let us sketch some properties of the lattice Heisenberg VOA $V_{\mathbb{Z}}$ and,  isomorphic to it, the fermionic superalgebra $\Lambda$.

On the VOA level, the pure algebraic construction of ${V_{\mathbb{Z}}}\cong{V_{\mathbb{Z}}}(P_0)$ may start with ${\Pi_0}\cong{\Pi(P_0)}\cong{\mathbb{C}[b_{-1},b_{-2},\ldots]}|0>$ and proceeds as follows ~\cite{EFDB01}:  For each $n\in{\mathbb{Z}}$ we introduce the formal space $\Pi_n$ which is an underlying space of an irreducible representation of the Heisenberg algebra $\mathcal{H}$ and which is given as ${\Pi_n}={\mathbb{C}}[b_{-1},b_{-2},\ldots]|n>$ . This means that the formal vectors $|n>$ have properties $b_k|n>=0$ for $k>0$ and  $b_0|n>=n|n>$. Then we introduce the space $V_{\mathbb{Z}}={\oplus}_{n\in{\mathbb{Z}}}\Pi_n$ and operators $\phi_k,\phi^{*}_k$ , $k\in{\mathbb{Z}}$ which satisfy  $\phi_k|0>=0$ for $k\geq{0}$ and $\phi_k^{*}|0>=0$ for $k>0$  and which have the   following anticommutative relations:
\begin{equation}
[\phi_k,\phi_l]_{+}=0,\quad [\phi_k^{*},\phi_l^{*}]_{+}=0, \quad [\phi_k,\phi_l^{*}]_{+}={\delta}_{k,-l}
\end{equation} 

By means of these operators  we may obtain any vector $|n>$  from the vector $|0>$  as well as we can construct any Heisenberg generator  $b_n$. Namely, we have
\begin {equation}
|n>=\phi_{n-1}^{*}\ldots\phi_{1}^{*}\phi_{0}^{*}|0> \quad and \quad |-n>=\phi_{-n}\ldots\phi_{-1}|0> \quad for \quad n>0  
\end{equation}

as well as 
\begin{equation}
b_{n}={\sum}_{k\in{\mathbb{Z}}}\phi_k^{*}\phi_{n-k} \quad for \quad n\neq{0}, \qquad b_{0}={\sum}_{k=0}^{\infty}\phi_{-k}^{*}\phi_{k}-{\sum}_{k<0}\phi_{k}\phi_{-k}^{*} 
\end{equation}

Similarly as all informations about the vertex algebra ${\Pi_{0}}\cong{\Pi(P_0)}$, $\Pi_{0}\subset{V_{\mathbb{Z}}}$, were contained in the vertex operator   $Y(b_{-1},z)=b(z)={\sum}_{n\in{\mathbb{Z}}}b_{n}z^{-n-1}$ now we have that all informations about ${V_{\mathbb{Z}}(P_0)}\cong{V_{\mathbb{Z}}}$ are contained in the fields
\begin{equation}
V_1(z):=Y(|1>,z)={\sum}\phi_n^{*}z^{-n} \qquad and \quad V_{-1}(z):={\sum}\phi_{n}z^{-n-1}
\end{equation}

The Heisenberg vertex operator $b(z)$ can be express in terms of the above vertex operetors  as their normally ordered product i.e.
\begin{equation}
b(z)=:V_1(z)V_{-1}(z):
\end{equation}

   Let us return to our Riemann surface $\Sigma$ and to a point $P_0$ on it. The appearance of the Heisenberg vertex operator algebra $\Pi(P_0)$ was naturally related to the algebraic structure of $\mathbb{C}((z))\cong{\mathcal{K}_{P_0}}$ equipped with the non-degenerate antisymmetric scalar product on it.  Since we may identify $\mathcal{K}_{P_0}$ with the space $\Gamma(\Delta_{\epsilon}^{*}(P_0),{\mathcal{O}(\underline{\mathbb{C}})})$ of holomorphic sections of the holomorphically trivial line bundle $\underline{\mathbb{C}}\in{Picc_0}$ restricted to ${U_{0}\cap{U_1}}={\Delta^{*}_{\epsilon}(P_0)}$  we may ask a question whether there exists a holomorphic line bundle $\xi\in{Picc_n}$, $n>0$ whose sections over $\Delta^{*}_{\epsilon}(P_0)$ carry an interesting algebraic structure which naturally leads to some vertex operator algebra as in the Heisenberg case.  We will see that the lattice algebra $V_{Z}$ (more exactly, the isomorphic to it, fermionic superalgebra $\Lambda$) is  naturally associated to $Picc_{g-1}\subset{Picc(\Sigma)}$ as well as that this correspondence provides an important tool for characterization  of any concrete compact Riemann surface. 

 Before we will introduce $\Lambda$ let us recall what happens in a finite dimentional case.  When we have a finite 2k-dimentional vector space which is equipped with a nondegenerate symmetric bilinear form we may introduce the Clifford product and we may construct the associated  Clifford algebra $\textsl{CL}(V)$. Its underlying vector space is spanned by all multivectors with respect to the Clifford multiplication (and hence it is isomorphic to the underlying space of the exterior algebra of $V$). The underlying space of an irreducible representation for $\textsl{Cl}(V)$  is given by a space of the algebraic spinors i.e. by a minimal left ideal $\textsl{Cl}(V)f$ of $\textsl{Cl}(V)$  which may be determined by an isotropic k-vector, say,  $f=e_{1},\ldots,e_k$.  Each element of V itself  can be built out of  algebraic spinors in a quite natural way ~\cite{CLC97}.

Let us go back to the point $P_{0}\in{\Sigma}$  and to the mutually dual spaces ${\mathcal{K}_{P_0}}\cong{\mathbb{C}((z))}$ and ${\Omega_{\mathcal{K}}}\cong{\mathbb{C}((z))dz}$. We have seen that the dual pairing given by the residue defines antisymmetric product $(f,g)$ on $\mathbb{C}((z))$ which further, when we work with meromorphic pictures P-II or P-III, leads to the Heisenberg Lie algebra $\mathcal{H}$ and then to the Heisenberg VOA  $\Pi(P_0)$. Now let us introduce another scalar product $\left\langle ,\right\rangle$ on $\mathbb{C}((z))$  which we will define as 
\begin{equation}
{\left\langle f,g\right\rangle}=Res_{z=0}fgdz
\end{equation} 
 
This product is symmetric (hence it does not define any extension of $\mathbb{C}((z))$) and it decomposes the whole space into two null subspaces, one of which is given by all principal parts.  We will denote the space $\mathbb{C}((z))$ that is equipped with the symmetric product $\left\langle ,\right\rangle$ by $H$ and its  null subspaces by $H_{+}$ and $H_{-}$ respectively. More precisely we consider
\begin{equation}
H={{\lim}_{{\epsilon}\rightarrow{0}}{\mathcal{O}}(\Delta_{\epsilon}^{*})}\cong{{\mathbb{C}}[[z]]}\oplus{{z^{-1}\mathbb{C}}[[z^{-1}]]}\equiv{H_{+}\oplus{H_{-}}}
\end{equation}

Similarly as in a finite-dimentional case, the symmetric product in H allows us to introduce the Clifford product and to construct the Clifford algebra $\textsl{Cl}(H)$. Let us denote the basis vectors $z^n$ of $H$ as $\nu_{n}:={z^n}\in{H}$.  The subspaces $H_{+}$ and $H_{-}$ are maximally isotropic (with respect to $\left\langle ,\right\rangle$) subspaces and a maximal totally isotropic multivector generates a minimal left ideal of $\textsl{Cl}(H)$. We will take   ${\Psi_0}=\nu_{-1}\nu_{-2}\nu_{-3}\ldots$ as such totally isotropic  multivector and we will denote it (similarly as before) as the element ${|0>}\in{\textsl{Cl}(H)}$. Since the Clifford product involves two products i.e.
\begin{equation}
\nu_{m}\nu_{n}=\frac{1}{2}\left\langle {\nu}_m,{\nu}_{n}\right\rangle+{\nu}_{m}\wedge{\nu}_n
\end{equation}

(where $\wedge$ denotes the wedge product  i.e. $\wedge:{H\times{H}}\rightarrow{H{\wedge}H}\subset{\Lambda}^{\infty}H$) it is natural to introduce (as in a finite dimentional case) two homomorphisms ${\widetilde{\rho}},{\widetilde{\rho^{*}}}:H\rightarrow{End{\textsl{Cl}(H)}}$ that are induced by these two   products. Thus we have
\begin{equation}
{\widetilde{\rho}({\nu}_n){\nu}_k}={\nu}_n{\wedge}{\nu}_k \quad and \quad {\widetilde{\rho^{*}}({\nu}_n)}{\nu_k}=\left\langle {\nu}_n,{\nu}_k\right\rangle
\end{equation}

We will denote these endomorphisms as follows
\begin{equation}
{\psi}_n^{*}:={\widetilde{\rho}({\nu}_{-n})} \qquad and \quad {\psi}_{n}:={\widetilde{\rho^{*}}(\nu_{-n-1})}
\end{equation}

We see immediately that these endomorphisms satisfy the following anticommutative relations
\begin{equation}
[\psi_n,\psi_m]_{+}=0, \quad [\psi_n^{*},\psi_m^{*}]_{+}=0, \qquad [\psi_n,\psi_m^{*}]_{+}={\delta}_{n,-m}
\end{equation}

and hence themselves generate   a Clifford algebra which is isomorphic to (and may be identify with) the Clifford algebra $\textsl{Cl}(H)$. We will denote this Clifford algebra by $\textsl{Cl}$. Since we have
\begin{equation*}
{\psi}_n^{*}{\nu_{-1}\nu_{-2}\ldots}\equiv{\psi_n^{*}|0>}=0, \quad for \quad n>0,\quad and \quad {\psi_n}|0>=0 \quad for \quad n\geq{0}
\end{equation*}

we obtain the, so called, fermionic Fock representation $\Lambda$ of $\textsl{Cl}$. It is generated by $\Psi_{0}=|0>={\nu_{-1}\nu_{-2}\ldots}$ and it has a basis that consists of vectors: 
\begin{equation}
\psi_{n_1}\ldots{\psi}_{n_k}\psi_{l_1}^{*}\ldots{\psi}_{l_m}^{*}|0> \quad with\quad {n_1}<{n_2}<...<{n_k}<0, \quad {l_1}<{l_2}<...<{l_m}\leq{0} 
\end{equation}

determined by the appropriate monomials in $\textsl{Cl}$. By setting the parity of these monomials to be equal to $(k+m)mod{2}$ we define the superspace structure on $\Lambda$. Now, similarly as before, we introduce a VOA structure  on $\Lambda$ by introducing the fields
\begin{equation}
{\Psi(z)}={\sum}\psi_{n}z^{-n-1} \qquad and \quad {\Psi}^{*}(z)={\sum}{\psi_n}^{*}z^{-n}
\end{equation}

These fields  carry all informations about the vertex operator superalgebra $\Lambda$. In particular, the normally ordered product $h(z)$ of these field, i.e. 
\begin{equation}
h(z)=:{\Psi}^{*}(z){\Psi}(z): \quad is \quad h(z)={\sum}{h_N}{z^{-N-1}}
\end{equation}

has  the coefficients ${h_N}={\sum}{\psi}_k^{*}\psi_{N-k}$ for $N\neq{0}$ and ${h_0}={\sum}_{k=1}^{\infty}{\psi^{*}_{-k}{\psi_k}}+{\sum}_{k\leq{0}}{\psi_k}{\psi^{*}_{-k}}$.  It is easy to check that these coefficients ${h_n}, n\in{\mathbb{Z}}$  satisfy the Heisenberg commutation relations analogous to $(3.5)$ (i.e.$[h_n,h_k]=n\delta_{n,-k}$ and   $h_0{\Psi_m}=m{\Psi_m}$)

We notice immediately the similarity between the formulas $(4.1)-(4.5)$ valid in  $V_{\mathbb{Z}}$ and the formulas $(4.11)-(4.14)$ . In fact, the anticommutative relations  $(4.1)$ for the operators ${\phi_n},{\phi}_n^{*}$ imply that these operators also generate a Clifford algebra and this Clifford algebra is isomorphic to $\textsl{Cl}$. In other words , we obtain a natural isomorphism $\sigma$ between the irreducible representations for these two Clifford algebras  which for any $n>0$  gives us
\begin{equation}
{\Psi}_{-n}:={\psi_{-n}\psi_{-n+1}...\psi_{-1}|0>}\stackrel{\sigma}{\rightarrow}{|-n>}\in{V_{\mathbb{Z}}} 
\end{equation}

Thus, in particular, we have ${\psi}_{-1}|0>\stackrel{\sigma}{\rightarrow}{\phi}_{-1}|0>=|-1>$,  etc. Similarly
\begin{equation}
{\Psi}_{n}:={\psi_{-n+1}^{*}...\psi_{-1}^{*}\psi_{0}^{*}|0>}\stackrel{\sigma}{\rightarrow}{|n>}\in{V_{\mathbb{Z}}}  
\end{equation}

 (in particular, ${\psi_0^{*}}\stackrel{\sigma}{\rightarrow}{\phi_0^{*}}|0>=|1>$).  Moreover we have 
\begin{equation}
{\Psi(z)}\stackrel{\sigma}{\rightarrow}{V_{-1}(z)} \qquad and \quad {\Psi^{*}(z)}\stackrel{\sigma}{\rightarrow}{V_1(z)}
\end{equation}

 Hence the decomposition $V_{\mathbb{Z}}={\oplus}_{n\in{\mathbb{Z}}}{\Pi_n}$ into subspaces $\Pi_{n}={\mathbb{C}[b_{-1},b_{-2},...]}|n>$  of irreducible representations of the Heisenberg Lie algebra $\mathcal{H}$ corresponds to the analogous decomposition ${\Lambda}={\oplus}_{n\in{\mathbb{Z}}}{\Lambda}^{(n)}$.  Each $\Lambda^{(n)}$ forms an irreducible representation of the oscillator algebra (generated by $h_n$)  which is based on the vector $\Psi_n$ instead of on $\Psi_{0}=|0>$.  Usually we view  the isomorphism $\sigma:{\Lambda}\rightarrow{V_{\mathbb{Z}}}$ as the direct sum of maps ~\cite{KAR87}
 \begin{equation}
 {\sigma}={\oplus}{\sigma_m}; \qquad {\sigma_m}:{\Lambda}^{(m)}\rightarrow{\Pi_m} 
 \end{equation}
 
 \subsection{HLB and vertex operator superalgebras $\Lambda$ and $V_{\mathbb{Z}}$} 
   We may notice that the formal product ${\left\langle f,g\right\rangle}=Rez_{z=0}f(z)g(z)dz$ in $\mathbb{C}((z))$ introduced earlier appears quite naturally when we  identify the space $H={\oplus}_{n\in{\mathbb{Z}}}{\mathbb{C}\nu_n}$ with the space of half-forms attached at ${P_0}\in{\Sigma}$  i.e. when we  view each $\nu_k$  as 
\begin{equation*}
{\nu_k}={z^k}dz^{\frac{1}{2}}\in{H}
\end{equation*} 

However to do this, first we must  fix some non-singular even spinor bundle ${\xi_{\beta}}\in{Picc_{g-1}}$  and any  local trivialization of the bunndle $\xi_{\beta}$ over ${U_0}={\Delta_{\epsilon}(P_0)}$. Any such choice identifies the space $H$ with the space of half forms at $P_0$. In other words, similarly as in the previous case, the essential   identifications 
\begin{equation}
{\mathbb{C}((z))}\cong{\Gamma({\Delta_{\epsilon}^{*}}(P_0),{\mathcal{O}(\underline{\mathbb{C}})})}, \qquad {\mathbb{C}((z))dz}\cong{\Gamma(\Delta_{\epsilon}^{*}(P_0),\mathcal{O}^{1,0})}
\end{equation}

 lead  to the Heisenberg vertex operator algebra $\Pi(P_0)$,  the   fixing of a non-singular even spinor line bundle $\xi_{\beta}\in{Picc_{g-1}}$,  gives us the identification
\begin{equation}
H={\oplus}_{n\in{\mathbb{Z}}}{\mathbb{C}\nu_{n}}\cong{\Gamma({\Delta}_{\epsilon}^{*}(P_0),{\mathcal{O}}(\xi_{\beta}))}
\end{equation}

which produces the vertex algebra $\Lambda(P_0)$. A non-singular even spinor bundle over any compact Rieman surface always exists. It corresponds to a class of divisors in $\textit{Div}\Sigma\diagup{\textit{Div}{\mathcal{P}}}$ that is determined by the divisor 
\begin{equation*}
P_0^{-1}P_{1}...P_g \qquad with \quad {P_0}\neq{P_i},\quad i=1,...g
\end{equation*}

of a unique meromorphic section $\sigma_{\beta}$ of $\xi_{\beta}$ with a single simple  pole at $P_0$.   The integral divisor $D={{P_1}...P_g}$ is a nonspecial one. The space 
\begin{equation}
W:={\mathbb{C}[\Sigma\backslash{P_0}]}{\sigma_{\beta}}
\end{equation}

of all sections of the holomorphic line bundle $\xi_{\beta}\in{Picc_{g-1}}$ with unique pole at $P_0$ is analogous to the space ${\mathbb{C}[\Sigma\backslash{P_0}]}$ of all sections with unique pole at $P_0$ of the line bundle $\underline{\mathbb{C}}\in{Picc_0}$. (The obvious diffrence is that whereas the canonical line bundle and the bundle $\underline{\mathbb{C}}$ are naturally distinguished in the spaces $Picc_{2g-2}$ and $Picc_0$ respectively, we have to fix an even non-singular spinor bundle $\xi_{\beta}$ in $Picc_{g-1}$). 

 Since the identification of the space $H$ with the space of half forms at $P_0$ given by $(4.20)$ allows us to associate  the vertex superalgebra $\Lambda(P_0)\cong{\Lambda}$ to the point  $P_0$  we may  identify  the space ${\mathbb{C}[\Sigma\backslash{P_0}]}{\sigma}_{\beta}$  with some element $\tau\in{{\Lambda}^{(0)}}\subset{\Lambda}(P_0)$. 

To see this we introduce the infinite grassmannian $GrH$ i.e. the collection of closed subspaces $V\subset{H}$ such that the projection $pr_{-}:V\rightarrow{H_{-}}$ is a Fredholm of the index zero and $pr_{+}:V\rightarrow{H_{+}}$ is a Hilbert-Schmidt. The subspace $W$ of $H$ determined by $(4.21)$ forms a concrete element of $Gr{H}$ and when we use the Grassmann variables we obtain a natural correspondence between the space $W$ and some element ${\tau\in(\textsl{Cl}({H}))\Psi_0}\cong{\Lambda}(P_0)$. Morover, $\tau\in{\Lambda}^{(0)}$ and since operators $\psi_{n}={\widetilde{\rho^{*}}}(\nu_{-n-1})$ and $\psi^{*}_{m}={\widetilde{\rho}}(\nu_{-m})$ given by $(4.10)$  act naturally on $\tau$, the totality of the Pluecker relations in the Grassmann variables is equivalent to ~\cite{EFDB01}, ~\cite{VKA97}
\begin{equation}
Res_{z=0}({\Psi(z){\tau}}\otimes{\Psi^{*}(z){\tau}})=0 \quad in \quad {\Lambda(P_0)}
\end{equation}

On the other hand, when we use the fact that (by the isomorphism $\sigma$) we have ${\Lambda^{(0)}}\cong{\Pi_0}\subset{V_{\mathbb{Z}}}$  we may express the element ${\tau}\cong{W}\in{Gr{H}}$ in terms of the elements of $\Pi_0$. For this let us change the notation to a more traditional one ~\cite{KAR87}  and let us write  the representation space of the Heisenberg Lie algebra $\mathcal{H}$ as $\mathbb{C}[x_1,x_2,...]$ . Now the action of the operators ${b_n}\in{\mathcal{H}}$, $n\in{\mathbb{Z}}$ are given by ${b_n}={\frac{\partial}{\partial{x_n}}}$ for $n>0$ and the action of $b_{-n}$ is the multiplication by $``n{x}_{n}''$. Thus, when we write ${\tau}\cong{W}\in{GrH}$ in the variables $x_k$'s then the relation  $(4.22)$ becomes 
\begin{equation}
Res_{z=0}(e^{{\sum}_1^{\infty}{\frac{b_{-i}}{i}}k^i}e^{-{\sum}_1^{\infty}{\frac{b_j}{j}}k^{-j}}{\tau}(x))(e^{-{\sum}_1^{\infty}{\frac{b_{-i}}{i}}k^i}e^{{\sum}_1^{\infty}{\frac{b_j}{j}}k^{-j}}{\tau}(x'))=0
\end{equation}

 We have used the coordinate $k={\frac{1}{z}}$ on $\Delta_{\epsilon}^{*}(P_0)$ instead of $z$. It appears that we may express $\tau(x)\in{\mathbb{C}[x_1,x_2,...]}$ in terms of the $\Theta$-function as follows: Let us normalize the  (introduced earlier)  meromorphic differentials $\tau^{(n+1)}_{P_0}$,  $n\geq{1}$ with the   single pole of order $n+1$ at $P_0$ and with all a-periods zero  by the condition
\begin{equation}
{\int_{P_0}^{P}\tau^{n+1}_{P_0}}=z^{-n}-2{\sum}_{j=1}^{\infty}Q_{nj}{\frac{z^j}{j}} 
\end{equation}

Let Q denote the matrix obtained by the coefficients $Q_{nj}$ above  and let $Q(x)={\sum}Q_{ij}x_{i}x_{j}$. Now ~\cite{EPR96}, ~\cite{IKR07}
\begin{equation}
\tau(x)=e^{Q(x)}\Theta(\textbf{B}^{\infty}x+e)
\end{equation}

The matrix $\textbf{B}^{\infty}$  was introduced earlier. The point $e\in{Jac\Sigma}$ is uniquely determined by the point ${P_0}\in{\Sigma}$ and by the spin bundle $\xi_{\beta}\in{Picc_{g-1}}$ (more precisely by the integral divisor $D$ introduced above).

  Now, the trivial bundle $\underline{\mathbb{C}}$ (in the Heisenberg case) is  naturally determined whereas, in the latter case, the spin bundle $\xi_{\beta}$ has to be chosen additionally to the choices of $P_0$ and of a local coordinate $z$ at $P_0$.  Since any additional choice accounts to giving some additional information it is not surprising that fixing $\xi_{\beta}$  provides informations about isospectral deformations of a normalized differential operator $\textbf{P}$  whose spectal curve is given by our Riemann surface $\Sigma$.  In other words, both equations $(4.22)$ and $(4.23)$ supply solutions to a KP-hierarchy of soliton equation associated to the spectral curve $\Sigma$ and to the point ${P_0}\in{\Sigma}$.

To complete our presentation the description of $\Sigma$ as a spectral curve and its relation to a K-P hierarchy  is given in the appendix.

\section{APPENDIX}

\subsection{From operator to its spectral curve}
Let $\mathcal{D}$ denote the algebra of ordinary differential operators with analytic coefficients i.e. ${\mathcal{D}}=\{\sum_{k=0}^{N}a_k(x){(\frac{d}{dx})}^{k},N\geq{0}\}$.  Let ${\textbf{P}}\in{\mathcal{D}}$  be a monic operator of order $N>1$ in the normalized form. If the commutative ring $\mathcal{A}_{\textbf{P}}\subset{\mathcal{D}}$ of all differential operators $\textbf{Q}\in{\mathcal{D}}$ that commute with $\textbf{P}$ has rank 1 then the problem of multiplicities of the eigenvalues of the operator $\textbf{P}$ is completely resolved.  This means that the space $Spec\textbf{P}$ of all eigenvalues of $\textbf{P}$ forms (after resolving multiplicities) N-sheeted covering  of  $\mathbb{C}$ and hence its completion $\overline{Spec\textbf{P}}$ is a compact Riemann surface, say ${\Sigma}={Spec\textbf{P}}\cup{\{\infty\}}$. We call this surface the spectral curve of $\textbf{P}$. Now,  $Spec\textbf{P}$ is an object af analysis and it is more convinient in this case to pass to the  world of algebraic geometry. Namely, since by the Schur lemma the algebra $\mathcal{A}_{\textbf{P}}$ is  commutative  we may consider the space $Spec{\mathcal{A}_{\textbf{p}}}$ of all prime ideals of $\mathcal{A}_{\textbf{P}}$. There is a natural map ${Spec\textbf{P}}\rightarrow{Spec{\mathcal{A}}_{\textbf{P}}}$ which allows us to identify the surface $\Sigma$ with ${\Sigma}\cong{\overline{Spec{\mathcal{A}_{\textbf{P}}}}}={Spec{\mathcal{A}_{\textbf{P}}}}\cup{\{\infty\}}$ . The point  ${P_0}\in{\Sigma}$ corresponds to $\{\infty\}$  and hence the space of all holomorphic functions on an algebraic variety $Spec{\mathcal{A}_{\textbf{P}}}$ (which is given by the algebra $\mathcal{A}_{\textbf{P}}$ itself) must coincide with the space ${\mathbb{C}}[\Sigma{\backslash{P_0}}]$. Moreover, gluing together simultaneous eigenspaces (each of which is one dimentional) of ${\mathcal{A}_{\textbf{P}}}$ at each point of $\overline{Spec{\mathcal{A}_{\textbf{P}}}}\cong{\Sigma}$ produces a holomorphic line bundle ${\xi}\in{Picc_{g-1}}(\Sigma)$.

It occurs that we may deform the operator ${\textbf{P}}\equiv{\textbf{P}(0)}$ without changing the spectral curve $\Sigma$. Such deformation ${\textbf{P}}={\textbf{P}(0)}\rightarrow{\textbf{P}(t)}$ (where $t=(t_1,t_2,...)$ denotes the parameters of isospectral deformations) is associated only with a change of  holomorphic line bundles from  $\xi$ to some $\xi_{t}$ in $Picc_{g-1}$.

Any isospectral deformation ${\textbf{P}(t)}$ of the operator $\textbf{P}$ must satisfy the Lax equation
\begin{equation}
\frac{d{\textbf{P}}(t)}{dt}=[\textbf{Q}(t),\textbf{P}(t)]
\end{equation} 

 for some ${\textbf{Q}(t)}\in{\mathcal{D}}$.  It is convinient to introduce a Lie algebra $\mathfrak{E}$ of pseudo-differential operators with analytic (or meromorphic ) coefficients
 \begin{equation*}
 \mathfrak{E}=\{\sum_{k=N}^{-\infty}a_{k}(x)\partial^{k};N\in{\mathbb{Z}}, {\partial}={\frac{d}{dx}}\}
 \end{equation*}

 We will consider the  pseudodifferential operator ${\textbf{P}}^{\frac{1}{N}}$ which we will denote by $L=L(0)$. Now the condition for isospectral deformations given by $(5.1)$  may be written in the following equivalent form :
 \begin{equation}
 \frac{dL(t)}{dt}=[\textbf{Q}(t),L(t)]
 \end{equation}
 
 It was shown by Gelfand and Dikki that all possible isospectral deformations are given by the Lax equations
 \begin{equation}
 \frac{dL(t)}{dt_k}=[L^k(t)_{+},L(t)] ; \quad k=2,3,...
 \end{equation}
 
 where ${L^k(t)_{+}}\in{\mathcal{D}}$ denotes the differential operator part of $L^k(t)$ and $L^k(t)_{-}$ is its pure pseudo-differential part. Writing the pseudo-differential operator $L(t)$ explicitely in the form
 \begin{equation}
 L(t)= \partial+\sum_{i=1}^{\infty}u_i(t)\partial^{-i}, \quad \partial={\frac{d}{dx}}, \quad t=(t_1,t_2,...), \quad t_1=x
 \end{equation}
 
 we immediately obtain that $L(t)$ is a solution to $(5.3)$ (i.e.describes isospectral deformation $\textbf{P}(t)=L(t)^{N}$) if and only if the coefficients $u_i(t)$'s are solutions to a sequence of nonlinear partial differential equations which form the so called  K-P hierarchy. ~\cite{IKR07}.  The nonlinear equations of K-P hierarchy are equivalent to some set of linear differential equations. To see this we introduce a pseudo-differential operator  $S(t)$ by the following condition: ${\partial}=S^{-1}(t)L(t)S(t)$.  The operator $S(t)$ has the form
 \begin{equation}
 S(t)=1+\sum_{j=1}^{\infty}s_j(t)\partial^{-j}
 \end{equation}
 
This operator is uniquely determined by $L$ (up to conjugation by elements $\widetilde{S}\in{\partial^{-1}}{\mathbb{C}}[[\partial^{-1}]]$).  The solutions to the nonlinear equations of K-P hierarchy are equivalent to the solutions of the system
 \begin{equation}
 S(t)\partial{S(t)}^{-1}=L(t) \quad and \quad \frac{\partial{S(t)}}{\partial{t_k}}=-L(t)^k_{-}S(t)
 \end{equation}
 
 which can be rewriten in terms of $w(t,z):=S(t)exp\sum_{i=1}^{\infty}t_{i}z^{-i}$, ${t_1}=x$, as the following linear differential equations
 \begin{equation}
 \frac{\partial{w(t,z)}}{\partial{t_k}}=L^k(t)_{+}w(t,z)
 \end{equation}.
 
  The function  $w(t,z)$ is called a wave function ~\cite{TSH86}, ~\cite{BAD81} .  It must satisfies the eigenvalue  equation
  \begin{equation*}
  L(t)w(t,z)=z^{-1}w(t,z)
  \end{equation*}
  which means that $(\frac{1}{z})^N$ is the eigenvalue of the operator $\textbf{P}$ for $z$ in ${U_0}={\Delta_{\epsilon}}(P_0)$
  
    If we introduce the quotient space $\textbf{V}$ of $\mathfrak{E}$ by its maximal ideal $\mathfrak{E}x$ generated by the set $\{\textbf{P}x; \textbf{P}\in{\mathfrak{E}}\}$ then the pseudodifferential operator $S(0)$ determines a unique subspace $W\subset{\textbf{V}}$ which satisfies 
 \begin{equation}
 dim(W\cap{\partial^{-1}\mathbb{C}[[\partial^{-1}]]})=dim({\mathbb{C}[\partial]\diagup{W\cap{\mathbb{C}}[\partial]}})
 \end{equation}
 
 The subspace  W is almost equal to $\mathbb{C}[\partial]$ (the difference is finite dimentional) and hence $W\in{Gr{\textbf{V}}}$.  So, we have the following correspondences
 \begin{equation}
 Gr{\textbf{V}}\ni{W_t}\leftrightarrow{S(t)}\rightarrow{L(t)=S(t)\partial{S(t)}^{-1}}
 \end{equation} 
 
 where,( by the construction above) $S(t)\cong{w(t,z)}$.  Since  we have the standard relation between the wave function and the $\tau$-function   the Pluecker relations in the Grassmann variables for $W_{t}\in{Gr{\textbf{V}}}$ imply the Hirota bilinear equations for $\tau$. These Hirota equations  are equivalent to $(4.22)$ exactly.
 
 \subsection{From a surface $\Sigma$ to operators $\textbf{P}(t)$}
 Given a compact Riemann surface $\Sigma$. Let $\xi_{\beta}$ be an even non-singular spinor bundle in ${Picc_{g-1}}$ whose unique section $\sigma_{\beta}$ has divisor $P_0^{-1}D$ (with a  nonspecial, integral divisor $D={P_1...P_g}$). Let $u_{P_0}:\Sigma\rightarrow{Jac\Sigma}$ denote the Abel map and let $e={-u_{P_0}}(D)-\Delta$ where $\Delta$ is the vector of Riemann constants corresponding to $P_0$. These data determine, the so called Baker-Akhiezer function ~\cite{BAD81}
 \begin{equation}
 \Psi(t,P)=exp(\sum_{i=1}^{\infty}{t_i}\int_{P_0}^{P}{\tau_{P_0}^{i+1}})\frac{\Theta(\textbf{B}^{\infty}t+e-u_{P_0}(P))\Theta(e)}{\Theta(e-u_{P_0})\Theta(\textbf{B}^{\infty}t+e)}
 \end{equation}
 
 which, after introducing a local coordinate z on $\Delta_{\epsilon}(P_0)$ vanishing at $P_0$ can be rewritten in terms of the wave function $w(t,z)=e^{\sum{t_i}{z^{-i}}}(1+\sum_{j=1}^{\infty}w_j(t)z^j)$.  Now we may view the function $w(t,z)$ as the image of  some pseudo-differential operator $S(t)=1+{\sum_{i=1}^{\infty}w_i(t)\partial^{-1}}$ acting on $e^{\sum{t_i}z^{-i}}$, ${t_1}=x$,  i.e. 
 \begin{equation}
 w(t,z)=S(t)e^{\sum{t_i}z^{-i}}
 \end{equation} 
 
 Since this latter relation determines the operator $S(t)$ we  immediately uncover the world of pseudodifferential operators with   the operators  $L(t)=S(t)\partial{S(t)^{-1}}$, with   the local coordinate z on $U_0={\Delta_{\epsilon}}(P_0)$ providing the eigenvalues $k^N={(\frac{1}{z})^N}$ for the differential operators $\textbf{P}=L^N$ and with the space $W\subset{\textbf{V}}$ corresponding exactly to the subspace $W\in{Gr}H$ given by $(4.21)$.  It is not difficult to check that $w(t,z)$ satisfies the equations $(5.7)$ and hence the soliton equations of the K-P hierarchy for the components $u_j(t)$  of $L(t)$ (given by $(5.4)$) have solutions that are expressed in terms of the $\Theta$-function associated to $\Sigma$ ~\cite{EACC88}
 
  Now, for any fixed $t=(t_1,t_2,..)$ the quotient of theta functions in $(5.10)$ is a multiplicative  multivalued function on $\Sigma$ whose divisor has degree zero. This means that it determines a unique holomorphic line bundle, say ${\Phi_t}\in{Picc_0}$.  The line bundles $\xi_{t_{k}}=\Phi_{t_k}\otimes{\xi_{\beta}}\in{Picc_{g-1}}$ correspond to a linear flow on $Jac\Sigma$ along the constant field determined by the k-th column of the matrix $\textbf{B}^{\infty}$.  It is clear that the space $\textbf{T}_{L}$ of effective parameters for deformations  $L(t)$ of $L$ (and hence for isospectral deformations $\textbf{P}(t)$ of $\textbf{P}$) is $\textbf{T}_{L}\cong{Jac\Sigma}$.

\end{document}